\documentclass[11pt]{amsart}

 \newtheorem{thm}{Theorem}[section]
 \newtheorem{cor}[thm]{Corollary}
 
 \newtheorem{prop}[thm]{Proposition}
 \theoremstyle{definition}
 \newtheorem{defn}[thm]{Definition}
 \theoremstyle{remark}
 
 \numberwithin{equation}{section}

 \newcommand{\norm}[1]{\left\Vert#1\right\Vert}
 
 \newcommand{\C}{\mathbb{C}}

\begin{document}

\title[]
 {Algebraic Elements and Invariant
Subspaces}

\author{Yun-Su Kim.  }


\address{Department of Mathematics, University of Toledo, Toledo,
Ohio, U.S.A. } \email{Yun-Su.Kim@utoledo.edu}


\dedicatory{}



\newpage

\begin{abstract}
 We prove that if a completely non-unitary contraction $T$ in $L(H)$ has a
non-trivial algebraic element $h$, then $T$ has a non-trivial
invariant subspace.

\end{abstract}

\maketitle

\section*{Introduction}
One of the most interesting open problems is the invariant
subspace problem. The invariant subspace problem is the question
whether the following statement is true or not: \vskip0.2cm

 Every bounded linear operator $T$ on a
separable Hilbert space $H$ of dimension $\geq{2}$ over $\C$ has a
non-trivial invariant subspace. \vskip0.2cm


We know that the invariant subspace problem is solved for all
finite dimensional complex vector spaces of dimension at least 2.
Thus, in this note, $H$ denotes a separable Hilbert space whose
dimension is infinite. Since it is enough to consider a
contraction $T$, i.e., $\norm{T}\leq{1}$ on $H$, in this note, $T$
denotes a contraction.

First, we focus on completely non-unitary contractions to use a
property of a multiplicity-free operator of class $C_0$, and we
consider
 \emph{algebraic elements with respect to} a completely non-unitary contraction $T$ introduced in \cite{K}.

If $T$ is a contraction, then\vskip0.2cm

(Case 1) $T$ is a completely non-unitary contraction with a
\emph{non-trivial algebraic element}, or

(Case 2) $T$ is a completely non-unitary contraction without a
non-trivial algebraic element, that is, every non-zero element in
$H$ is \emph{transcendental with respect to }$T$, or

(Case 3) $T$ is not completely non-unitary.\vskip0.2cm

In this note, we discuss the invariant subspace problem for
operators of (Case 1) or (Case 3).
 By using a classification of the invariant subspaces of a
multiplicity-free operator of class $C_0$ (\cite{B1}), in Theorem
\ref{17}, we prove that, for a completely non-unitary contraction
$T$, if $T$ has a non-trivial algebraic element $h$, then $T$ has
a non-trivial invariant subspace. By Theorem \ref{17}, we conclude
that every $C_{0}$-operator has a non-trivial invariant subspace
(Corollary \ref{14}). It follows that every nilpotent operator has
a non-trivial invariant subspace (Corollary \ref{15}).


In Corollary \ref{10}, we prove that if $T_{1}$ is not a
completely non-unitary contraction, then $T_{1}$ also has a
non-trivial invariant subspace.

Thus, we answer to the invariant subspace problem for the (Case 1)
and (Case 3) in Theorem \ref{17} and Corollary \ref{10}
respectively. Therefore, to answer to the invariant subspace
problem, it suffices to answer for (Case 2).

We do not consider operators of (Case 2) in this note, and leave
as a question;

\vskip0.2cm

\textbf{Question}.  Suppose that every non-zero element in $H$ is
transcendental with respect to $T$. Then, does the operator $T$
have a non-trivial invariant subspace?



\section{Preliminaries and Notation}\label{13}

In this note, $\C$, $\overline{M}$ and $L(H)$ denote the set of
complex numbers, the (norm) closure of a set $M$, and the set of
bounded linear operators from $H$ to $H$ where $H$ is a separable
Hilbert space whose dimension is infinite, respectively.

If $T\in{L(H)}$ and $M$ is an invariant subspace for $T$, then
$T|M$ is used to denote the restriction of $T$ to $M$.

\subsection{A Functional Calculus.}\label{11}

 Let $H^{\infty}$
be the Banach space of all (complex-valued) bounded analytic
functions on the open unit disk $\textbf{D}$ with supremum norm
\cite{S2}. A contraction $T$ in $L(H)$ is said to be
\emph{completely non-unitary } provided its restriction to any
non-zero reducing subspace is never unitary.


B. Sz.-Nagy and C. Foias introduced an important functional
calculus for completely non-unitary contractions.
\begin{prop}\label{12}Let $T\in{L(H)}$ be a completely non-unitary
contraction. Then there is a unique algebra representation
$\Phi_{T}$ from $H^{\infty}$ into $L(H)$ such that :\vskip0.2cm

(i) $\Phi_{T}(1)=I_{H}$, where $I_{H}\in{L(H)}$ is the identity
operator;

(ii) $\Phi_{T}(g)=T$, if $g(z)=z$ for all $z\in\textbf{D}$;

(iii) $\Phi_{T}$ is continuous when $H^{\infty}$ and $L(H)$ are
given the weak$^\ast$-

\quad\quad topology.

(iv) $\Phi_{T}$ is contractive, i.e.
$\norm{\Phi_{T}(u)}\leq\norm{u}$ for all
$u\in{H^{\infty}}$.\end{prop}

We simply denote by $u(T)$ the operator $\Phi_{T}(u)$.

B. Sz.- Nagy and C. Foias \cite{S2} defined the \emph{class
$C_{0}$} relative to the open unit disk \textbf{D} consisting of
completely non-unitary contractions $T$ on $H$ such that the
kernel of $\Phi_{T}$ is not trivial.  If $T\in{L(H)}$ is an
operator of class $C_{0}$, then \begin{center}$\ker
\Phi_{T}=\{u\in{H^{\infty}}:u(T)=0\}$\end{center} is a
weak$^{*}$-closed ideal of $H^{\infty}$, and hence there is an
inner function generating ker $\Phi_{T}$. The \emph{minimal
function} $m_{T}$ of an operator $T$ of class $C_{0}$ is the
generator of ker $\Phi_{T}$, that is,
\begin{equation}\label{2}\ker \Phi_{T}=m_{T}H^{\infty}.\end{equation}

\subsection{Algebraic Elements}
To provide a sufficient condition for (non-trivial) invariant
subspaces, we will use the notion of \emph{algebraic elements} for
a completely non-unitary contraction $T$ in $L(H)$.

\begin{defn}\cite{K}
Let $T\in{L(H)}$ be a completely non-unitary contraction. An
element $h$ of $H$ is said to be \emph{algebraic with respect to}
$T$ provided that $\theta(T)h=0$ for some
$\theta\in{H^{\infty}}\setminus\{0\}$. If $h\neq{0}$, then $h$ is
said to be a \emph{non-trivial} \emph{algebraic element} with
respect to $T$.

If $h$ is not algebraic with respect to $T$, then $h$ is said to
be \emph{transcendental with respect to} $T$.

\end{defn}

   If $B$ is a closed subspace of $H$
generated by $\{b_{i}\in{H}:i=1,2,3,\cdot\cdot\cdot\}$, then $B$
will be denoted by $\bigvee_{n=1}^{\infty}b_{i}.$

\subsection{Blaschke Products}

For each $\alpha\in{\textbf{D}}$, the \emph{Blaschke factor}
$b_{\alpha}$ is a function defined by
\[b_{\alpha}(z)=\frac{{|\alpha|}}{\alpha}\frac{\alpha-z}{1-\overline{\alpha}z},\emph{ }z\in{\textbf{D}},\]
if $\alpha\neq{0}$, and $b_{0}(z)=z$ for $z\in{\textbf{D}}$. We
recall that a \emph{Blaschke Product} is a function of the form
\[b(z)=\prod_{j}b_{\alpha_{j}}(z),\emph{ }z\in{\textbf{D}},\]
where $\{\alpha_{j}\}_{j=0}^{\infty}$ is a sequence in
$\textbf{D}$ such that $\sum_{j}(1-|\alpha_{j}|)<\infty$.

 A function
$\mu:\textbf{D}\rightarrow\{0,1,2,\cdot\cdot\cdot\}$ is said to be
a \emph{Blaschke function} if
\[\sum_{\alpha\in{\textbf{D}}}\mu(\alpha)(1-|\alpha|)<\infty.\]
If $b$ is a Blaschke Product, then we have a Blaschke function
$\mu$ such that
\[b(z)=\prod_{\alpha\in{\textbf{D}}}(b_{\alpha}(z))^{\mu(\alpha)},\]
where $\mu(\alpha)$ represents the multiplicity of $\alpha$ as a
zero of $b$, and this Blaschke Product $b$ will be denoted by
$b_{\mu}$. Recall that a \emph{singular inner function} is
determined by a positive finite measure $\nu$ on
$\partial{\textbf{D}}$, singular with respect to Lebesgue measure,
via the formula
\[ s_{\nu}(\lambda)=\exp(-\int_{\partial{\textbf{D}}}\frac{\xi+\lambda}{\xi-\lambda}d\nu(\xi)),\]
for $\lambda\in{{\textbf{D}}}$ (\cite{B1}).

If $\theta$ is an inner function, then there exist a Blaschke
product $b$, a singular inner function $s$, and a constant
$\gamma$, $|\gamma|=1$, such that
\[\theta=\gamma{b}s.\]

Let us recall that the set of positive finite measures on
$\partial{\textbf{D}}$ has a lattice structure with respect to the
following relation;
$\nu\leq\nu^{\prime}\Leftrightarrow{\nu(A)\leq\nu^{\prime}(A)}$
for every Borel subset $A$ of $\partial{\textbf{D}}$, and the set
of Blaschke functions can also be organized as a lattice with
respect to the following relation;
$\mu\leq\mu^{\prime}\Leftrightarrow{\mu(z)\leq\mu^{\prime}(z)}$
for every $z\in{\textbf{D}}$ \cite{B1}.

Let $\theta$ and $\theta^{\prime}$ be two functions in
$H^{\infty}$. We say that $\theta$ \emph{divides}
$\theta^{\prime}$(or $\theta$$\mid$$\theta^{\prime}$) if
$\theta^{\prime}$ can be written as
$\theta^{\prime}=\theta\cdot\phi$ for some $\phi\in {H^{\infty}}$.
We will use the notation $\theta\equiv\theta^{\prime}$ if $\theta$
and $\theta^{\prime}$ are two inner functions that differ only by
a constant scalar factor of absolute value one. Thus, the
relations $\theta$$\mid$$\theta^{\prime}$ and
$\theta^{\prime}$$\mid$$\theta$ imply that
$\theta\equiv\theta^{\prime}$.

\begin{prop}(\cite{B1})\label{1}
Let $\mu$ and $\mu^{\prime}$ be Blaschke functions, $\nu$ and
$\nu^{\prime}$ singular measures on $\partial{\textbf{D}}$,
$\gamma$ and $\gamma^{\prime}$ complex numbers of absolute value
one, and set $\theta=\gamma{b_{\mu}}s_{\nu},$
$\theta^{\prime}=\gamma^{\prime}{b_{\mu^{\prime}}}s_{\nu^{\prime}}$.

Then, $\theta|\theta^{\prime}$ if and only if
$\mu\leq\mu^{\prime}$ and $\nu\leq\nu^{\prime}$.
\end{prop}

\section{The Main Results}

We recall that an operator $T$ ia said to be
\emph{multiplicity-free} if $T$ has a cyclic vector. The invariant
subspaces of a multiplicity-free operator of class $C_0$ have a
classification as following ;

\begin{prop}[\cite{B1}, Theorem 3.2.13]\label{7}
For every operator $T$ of class $C_0$, the following assertions
are equivalent;

(i) $T$ is multiplicity-free.

(ii) For every inner divisor $\phi$ of $m_{T}$ (that is, $\phi$ is
an inner function such that $\phi|m_{T}$), there exists a unique
invariant subspace $\emph{K}$ for $T$ such that
$m_{T|\emph{K}}\equiv\phi$.\vskip.1in

If $T$ is multiplicity-free, then the unique invariant subspace in
(ii) is given by $\emph{K}=\ker\phi(T)$.

\end{prop}

Note that $m_{T}$ always has two \emph{trivial inner divisors}
$\phi_{1}\equiv{1}$ and $\phi_{2}\equiv{m_{T}}$.

\begin{prop}\label{19}
Let $T:H\rightarrow{H}$ be a multiplicity-free operator of class
$C_0$ and $\phi$ be an inner divisor of $m_{T}$.

If $\phi$ is not a trivial inner divisor, then $\ker\phi(T)$ is a
non-trivial invariant subspace for $T$.

\end{prop}
\begin{proof}
Suppose $\ker\phi(T)=H$. Then, by equation (\ref{2}),
$m_{T}|\phi$. Since $\phi$ is an inner divisor of $m_{T}$,
$\phi\equiv{m_{T}}$ which is a contradiction.

Suppose that $\ker\phi(T)=\{0\}$. Since $\phi$ is an inner divisor
of $m_{T}$,
\begin{equation}\label{3}m_{T}=\phi\varphi,\end{equation} for a
function $\varphi\in{H^{\infty}}$.
 Since
$m_{T}(T)=\phi(T)\varphi(T)=0$ and $\ker\phi(T)=\{0\}$,
$\varphi(T)=0$. Thus, by equation (\ref{2}),
\begin{equation}\label{70}m_{T}|\varphi.\end{equation}
By equation (\ref{3}) and (\ref{70}), we have
$\varphi\equiv{m_{T}}$, and equation (\ref{3}) implies that
$\phi\equiv{1}$ which is a contradiction.

\end{proof}

\begin{thm}\label{17}
Let $T\in{L(H)}$ be a completely non-unitary contraction.

If $T$ has a non-trivial algebraic element $h$, then $T$ has a
non-trivial invariant subspace.
\end{thm}
\begin{proof}
Suppose that $T$ has a non-trivial algebraic element $h$. Then,
there is a non-zero function $\theta\in{H^{\infty}}$ such that
$\theta(T)h=0$.

Let $M=\bigvee_{n=0}^{\infty}T^{n}h,$ and $T_{1}=T|M$.

Since $\theta(T)(T^{n}h)=T^{n}(\theta(T)h)=0$ for any
$n=0,1,2,\cdot\cdot\cdot$, $\theta(T_{1})=0$. Thus, the operator
$T_{1}:M\rightarrow{M}$ is a multiplicity-free operator of class
$C_0$. If $m_{T_{1}}\equiv{1}$, then, by Proposition \ref{12},
\begin{equation}\label{20}m_{T_{1}}(T_{1})h=h.\end{equation} Since $m_{T_{1}}(T_{1})\equiv{0}$,
equation (\ref{20}) implies that $h=0$ which is a contradiction.
Thus, \begin{equation}\label{21}m_{T_{1}}\neq{1}.\end{equation}

By Proposition \ref{7}, for every inner divisor $\phi$ of
$m_{T_{1}}$, there exists a unique invariant subspace
$\emph{K}(=\ker\phi(T_{1}))$ for $T_{1}$.

By Proposition \ref{19}, if $\phi$ is not a trivial inner divisor
 of $m_{T_{1}}$, then $\ker\phi(T_{1})$ is a non-trivial invariant
subspace for $T_{1}$. Note that $\ker\phi(T_{1})$ is also a
non-trivial invariant subspace for $T$.

If $m_{T_{1}}$ has no non-trivial inner devisors, then, by
Proposition \ref{1} and (\ref{21}), $m_{T_{1}}$ must be a Blaschke
factor. Then, $m_{T_{1}}\equiv\frac{a-z}{1-\overline{a}z}$ for
$z\in\textbf{D}$ where $a$ is a complex number in $\textbf{D}$,
and in this case we clearly have $m_{T_{1}}|(z-a)$. Therefore,

\[(z-a)\in\ker \Phi_{T_{1}}=\{u\in{H^{\infty}}:u(T_{1})=0\}.\]
Thus, $(T-aI_{H})h=(T_{1}-aI_{H})h=0$ and so $h$ is an eigenvector
of $T$. Since $h\neq{0}$ and $\dim{H}\geq{2}$, the closed subspace
$M(=\bigvee_{n=0}^{\infty}T^{n}h)$ generated by the eigenvector
$h$ is a non-trivial invariant subspace for $T$.


\end{proof}
\begin{cor}\label{14}
Every $C_{0}$-operator has a non-trivial invariant subspace.

\end{cor}

\begin{cor}\label{15}
Every nilpotent operator has a non-trivial invariant subspace.

\end{cor}

Recall that an operator is \emph{subnormal} if it has a normal
extension. Thus, every normal operator is trivially subnormal.

\begin{prop}\cite{S3}\label{8}

Every subnormal operator has a non-trivial invariant subspace.

\end{prop}
Thus, we have the following well-known result;

\begin{cor}\label{10}
 If $T$ is not a completely non-unitary
contraction, then $T$ has a non-trivial invariant subspace.
\end{cor}
\begin{proof}
 Since $T$ is not a completely non-unitary
contraction, there is a reducing subspace $M(\neq\{0\})$ for $T$
such that $T|M$ is a unitary operator.

If $M$ is a non-trivial reducing subspace for $T$, then it is
done. If $M\equiv{H}$, then $T$ is unitary, and so $T$ is
subnormal. By Proposition \ref{8}, it is proven.

\end{proof}





------------------------------------------------------------------------

\bibliographystyle{amsplain}
\bibliography{xbib}
\end{document}